\documentclass[english]{IEEEtran}
\usepackage[T1]{fontenc}
\usepackage[latin9]{inputenc}
\usepackage{units}
\usepackage{amsmath}
\usepackage{graphicx}

\makeatletter

\providecommand{\tabularnewline}{\\}

\usepackage{cite}
\usepackage{fixltx2e}

\makeatother

\usepackage{babel}
\begin{document}

\title{Changes in Cascading Failure Risk with Generator Dispatch Method
and System Load Level}

\author{Pooya Rezaei, \emph{Student Member, IEEE}, Paul D. H. Hines, \emph{Member,
IEEE}%
\thanks{This work is supported in part by the US Dept.~of Energy award \#DE-OE0000447.

The authors are with the School of Engineering, University of Vermont,
Burlington, VT (e-mail: \{pooya.rezaei,paul.hines\}@uvm.edu).%
}}
\maketitle
\begin{abstract}
Industry reliability rules increasingly require utilities to study
and mitigate cascading failure risk in their system. Motivated by
this, this paper describes how cascading failure risk, in terms of
expected blackout size, varies with power system load level and pre-contingency
dispatch. We used Monte Carlo sampling of random branch outages to
generate contingencies, and a model of cascading failure to estimate
blackout sizes. The risk associated with different blackout sizes
was separately estimated in order to separate small, medium, and large
blackout risk. Results from $N-1$ secure models of the IEEE RTS case
and a 2383 bus case indicate that blackout risk does not always increase
with load level monotonically, particularly for large blackout risk.
The results also show that risk is highly dependent on the method
used for generator dispatch. Minimum cost methods of dispatch can
result in larger long distance power transfers, which can increase
cascading failure risk.\end{abstract}

\section{Introduction}

Cascading failure in power systems refers to a sequence of interdependent
outages that is initiated by one or more disturbances. Timely operator
intervention can often prevent a cascade from resulting in a large
blackout; however, large cascades occasionally occur and lead to major
blackouts, such as Aug.~2003 \cite{:2004} or Sept.~2011 \cite{:2012}.
Although large blackouts are low-probability events, they can have
catastrophic social and economic impacts. For this reason cascading
failure (CF, hereafter) risk assessment is increasingly required by
reliability regulations (e.g., from NERC \cite{NERC:2007}) and is
a focus of IEEE Power and Energy Society activities \cite{Vaiman:2012}.
State-of-the-art CF risk assessment methodologies are documented in
\cite{Vaiman:2012}.

In this paper, we used Monte Carlo simulation to estimate CF risk.
Monte Carlo methods are widely used for power system reliability evaluation
\cite{Allan:1996}. However, their application to the problem of estimating
CF risk is, on the other hand, much less established in the literature.
Standard reliability models typically only calculate the immediate
consequence of a sampled outage (such as the direct load shedding
that results from a generator outage). Estimating the additional risk
posed by the potential for cascading blackouts is more difficult for
several reasons. Firstly, the simulation of CF remains a difficult
problem, for which little validation data exist, and more research
is needed \cite{Vaiman:2012}. Secondly, even if a CF can be simulated,
the size (in terms of MW of load lost) of a CF can be at any scale,
which gives rise to the well documented power-law in blackout sizes
\cite{Carreras:2000,Carreras:2002}. Thirdly, the size of the search
space of all possible $n-k$ contingencies, where $n$ is the number
of components that might fail and $k$ is the number that did fail,
is enormous and grows exponentially with $n$ and $k$. Finally, the
combinations of outages (and operator errors) that typically trigger
CF are usually very low probability and very high impact, further
increasing the required computational effort. 

A few papers have used Monte Carlo sampling for CF risk estimation
\cite{Chen:2005,Kirschen:2004,Chen:2013,Kim:2013}, and some used
sampling techniques to reduce the computational cost of risk estimation.
The authors in \cite{Chen:2005} utilized Monte Carlo simulation and
then Importance Sampling to reduce the computational burden. Ref.~\cite{Kirschen:2004}
used correlated sampling in a Monte Carlo simulation to estimate the
stress level in the power system. The authors in \cite{Chen:2013}
showed that importance sampling together with variance reduction technique
can be used to the increase computational efficiency of CF risk estimation
by a factor of 2-4. In \cite{Kim:2013}, the Splitting method was
used to produce more substantial speed improvements for a system with
both continuous and discrete uncertainty. Non-sampling approaches,
such as branching process models \cite{Ren:2008,Dobson:2012} provide
efficient estimates of risk, but abstract away some details, such
as the ability to identify which transmission lines contribute to
risk estimates. The authors in \cite{Liao:2004} and \cite{Nedic:2006}
found a phase transition in CF risk when load level changes. 

In this paper, we are interested to study the impacts of pre-contingency
generator dispatch and load level on CF risk. Here, we use the expected
value of blackout sizes resulting from random contingencies as our
metric of risk. The rest of this paper is organized as follows. Section
\ref{sec:Monte-Carlo-Simulation} discusses the Monte Carlo simulation
method. Section \ref{sec:Pre-contingency-Dispatch} describes the
formulation and assumptions for the pre-contingency system. Section
\ref{sec:Simulation-and-Results} discusses the simulation and results.
Section \ref{sec:Monte-Carlo-Limitations} explains limitations of
the Monte Carlo approach and motivation for future research. Section
\ref{sec:Conclusions} provides our conclusions from this study.

\section{Monte Carlo Simulation\label{sec:Monte-Carlo-Simulation}}

In this paper, we estimate CF risk using Monte Carlo (MC) simulation.
In each MC iteration, we randomly choose a set of one or more transmission
line or transformer (branch) outages randomly, based on the failure
rate of each component. Ideally, one would select line outages using
a joint probability distribution function for line outages, accounting
for correlations among the outage probabilities. Since correlation
data are not generally available and correctly modeling correlated
line outages requires some care, we assume that line outages are independent
here. Accounting for correlations in line outage probabilities remains
for future work.

Given this assumption, the probability of two (or more) simultaneous
outages is the product of each line-failure probability. We use the
failure rate of transmission lines ($\lambda$, outages/year), and
assume that each failure lasts for 1 hour on average. Then, the probability
of a line failure at each iteration is computed by $p_{f,i}=\nicefrac{\lambda_{i}}{8760}$
for all lines, where 8760 is the number of hours in a year. 

Each random draw in the simulation produces a set of outaged lines
with a minimum size of zero. If the size of the outage is 2 or larger
(since the system is known to be initially $n-1$ secure), this contingency
is applied to our cascading failure simulator (CFS), which is explained
in \cite{Eppstein:2012} in detail. Following the standard MC approach,
we use the average (expected) blackout size in MW of load shedding
as our measure of risk. The expected value (average) is found by summing
over all event sizes and dividing by the number of MC iterations (including
events with zero blackout size and zero branch outages). 

In finding cascading failure risk, it is useful to separately consider
risk from events of different sizes. To do so, we add blackout sizes
within a certain size range in the numerator and divide by the number
of MC iterations, to find the risk associated with blackouts in different
size ranges.

\section{Pre-contingency Dispatch\label{sec:Pre-contingency-Dispatch}}

For the initial results in this paper, we computed the pre-contingency
power flow state for each load level using a Security-Constrained
DC Optimal Power Flow (SCDCOPF). As a result, each pre-contingency
network is $n-1$ secure for any single line outage. Although the
DC Optimal Power Flow (DCOPF) is a relatively simple linear programming
problem, a full SCDCOPF including all line outages as contingencies
can become computationally expensive, especially for larger systems,
because of extensive number of contingencies. To reduce the computational
effort, we solve a decomposed SCDCOPF based on the method proposed
in \cite{Li:2009}. The decomposed SCDCOPF is described as follows.

Initially, DCOPF is solved to find minimal generation cost dispatch
constrained by power flow equations and line flow limits, with the
following formulation:

\begin{eqnarray}
\min_{\mathbf{P}_{g},\mathbf{P}_{d}} &  & \boldsymbol{\mathbf{c}}_{g}^{T}\mathbf{P}_{g}-\boldsymbol{\mathbf{c}}_{d}^{T}\mathbf{P}_{d}\label{eq:obj}\\
\textrm{s.t.} &  & \mathbf{P_{\bar{r}}}=\mathbf{B_{\overline{rr}}\boldsymbol{\theta_{\bar{r}}}}\label{eq:dcpf}\\
 &  & \overset{\left|G\right|}{\underset{i=1}{\sum}}P_{g,i}=\overset{\left|D\right|}{\underset{j=1}{\sum}}P_{d,j}\label{eq:sigma}\\
 &  & \mathbf{F}=\mathbf{x}_{b}^{-1}\mathbf{A}\mathbf{\boldsymbol{\theta_{\bar{r}}}}\label{eq:flow}\\
 &  & -\mathbf{F}_{\max}\leq\mathbf{F}\leq\mathbf{F}_{\max}\label{eq:flow_lim}\\
 &  & 0\leq\mathbf{P}_{g}\leq\mathbf{P}_{g,\max}\label{eq:gen_lim}\\
 &  & 0\leq\mathbf{P}_{d}\leq\mathbf{P}_{d0}\label{eq:load_lim}
\end{eqnarray}
where $\mathbf{P}_{g}$ and $\mathbf{P}_{d}$ denote vectors of real
power generation and load in the network, and $G$ and $D$ are sets
of buses with generators and loads respectively. $\boldsymbol{\mathbf{c}}_{g}$
and $\boldsymbol{\mathbf{c}}_{d}$ are vectors of generator marginal
costs and the cost of load shedding at each bus, respectively (both
in $\nicefrac{\$}{MWh}$). $\mathbf{P_{\bar{r}}}$, $\mathbf{B_{\overline{rr}}\boldsymbol{}}$
and $\boldsymbol{\theta_{\bar{r}}}$ are respectively the vector of
real power injections, bus susceptance matrix, and the vector of bus
voltage angles for all buses except the reference bus. $\mathbf{F}$
denotes the line power flow vector. $\mathbf{x}_{b}$ is a matrix
with each diagonal entry representing the susceptance of each line,
and zero non-diagonal entries. $\mathbf{A}$ is the node-branch incidence
matrix, where the number of rows and columns are equal to the number
of branches and buses respectively. 

Constraints (\ref{eq:dcpf}) and (\ref{eq:sigma}) enforce the DC
power flow constraints, and constraint (\ref{eq:flow}) calculates
flows from bus voltage angles. Constraints (\ref{eq:flow_lim}), (\ref{eq:gen_lim})
and (\ref{eq:load_lim}) restrict line flows, real power generation
and load to be between their upper and lower bounds. Constraint (\ref{eq:load_lim}),
together with the second term in (\ref{eq:obj}), enables the possibility
of load shedding, which ensures that the problem is always feasible.
In order to ensure that load shedding does not occur unless absolutely
necessary, we set the entries of $\boldsymbol{\mathbf{c}}_{d}$ to
have large positive values, which are all greater than $\mathbf{c}_{g}$.
In this paper, we assume equal values of $c_{d}$ for all loads.

\begin{figure}
\begin{centering}
\includegraphics[width=0.5\columnwidth]{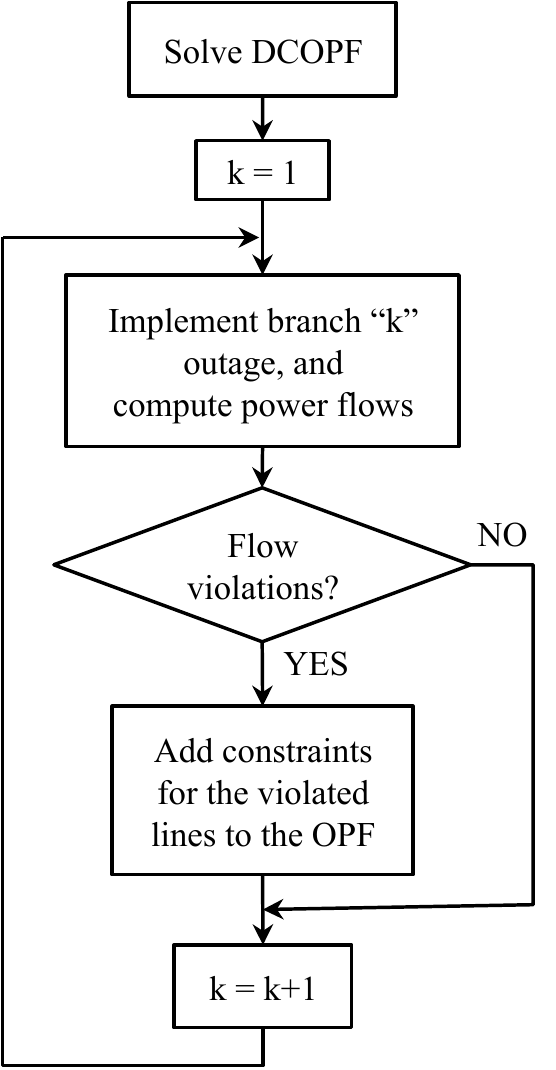}
\par\end{centering}

\caption{\label{fig:scopf_flowchart}One cycle of the decomposed SCDCOPF (based
on \cite{Li:2009})}

\end{figure}

In order to make each case $n-1$ secure, we add contingency constraints.
Here, we use the Line Outage Distribution Factors (LODF) matrix to
find post-contingency line flows after each line outage \cite{Wood:1996},
which is an $m\times m$ matrix, where $m$ represents the number
of branches. Assuming line $j$ is tripped in the network, each entry
$h_{ij}$ of the LODF matrix gives the relative change in flow on
line $i$ due to the outage of line $j$. Therefore, each post-contingency
flow constraint has the following form: 

\begin{equation}
-F_{i,\max}^{'}\leq f_{i}+h_{ij}f_{j}\leq F_{i,\min}^{'}
\end{equation}
where $F_{i,\max}^{'}$ denotes the short term rating of line $i$.
$f_{i}$ and $f_{j}$ denote the pre-contingency power flows on lines
$i$ and $j$ respectively. In order to solve a full SCDCOPF, one
can add as many as $m(m-1)$ contingency constraints to the problem.
However, explicitly adding these constraints makes problem prohibitively
computational expensive, especially for large $m$. To reduce the
computational cost, we implement a decomposed SCOPF, based on \cite{Li:2009},
in which contingency constraints are incrementally added to the problem
until the solution is $n-1$ secure. The flowchart of one cycle of
this algorithm is shown in Fig.~\ref{fig:scopf_flowchart}. Typically,
only 2 or 3 repetitions are needed to find an $n-1$ secure solution. 

\section{Simulation and Results\label{sec:Simulation-and-Results}}

\subsection{Test Networks\label{sub:Test-Networks}}

\begin{figure*}
\begin{centering}
\includegraphics[width=0.8\textwidth]{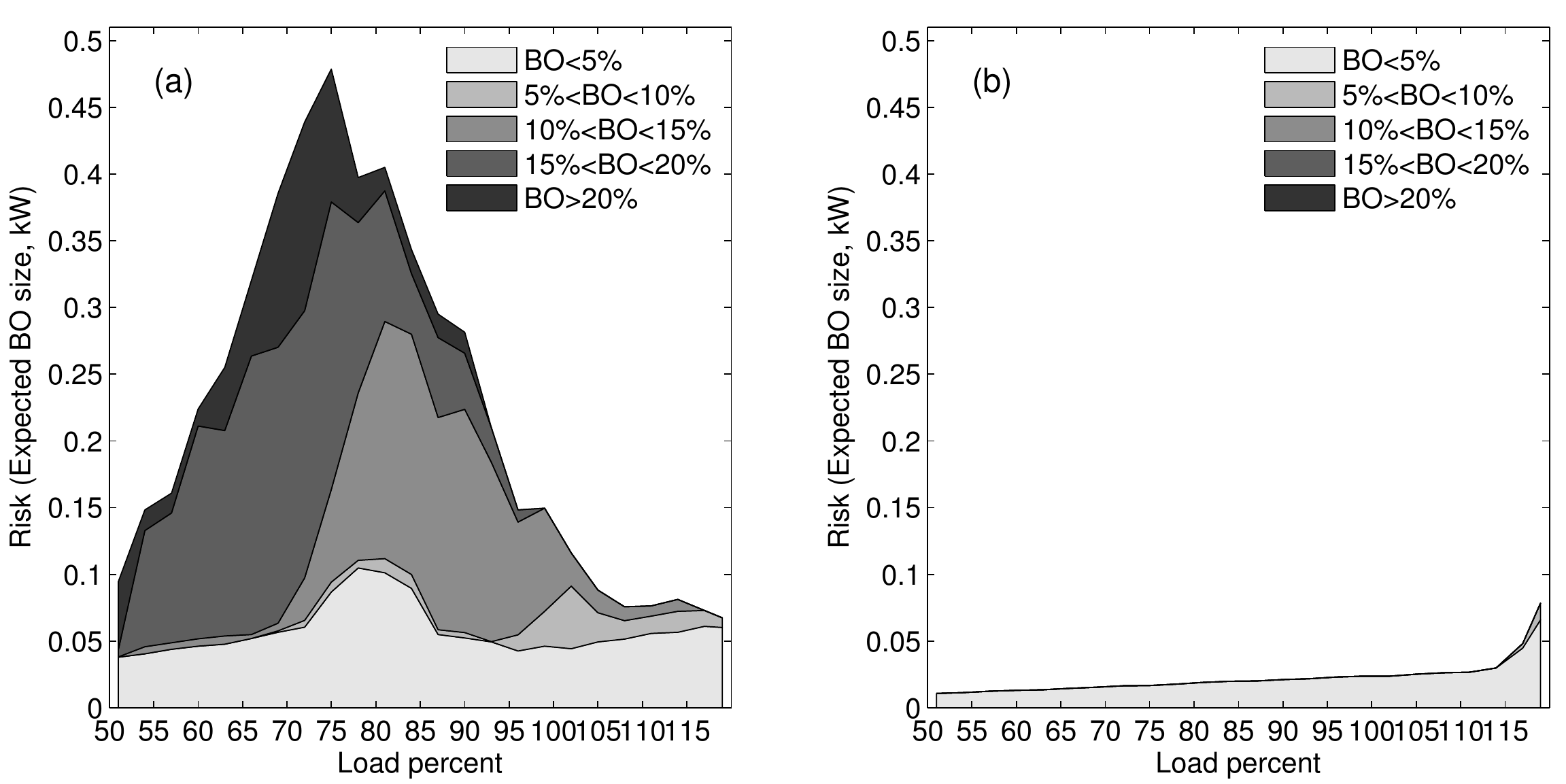}
\par\end{centering}

\caption{\label{fig:risk_RTS}CF risk (expected blackout size) for two pre-contingency
dispatch conditions: (a) SCDCOPF , and (b) Proportional}
\end{figure*}

We used two test cases to examine how cascading failure risk changes
with system load level. First, the 73-bus RTS-96, which has 120 branches
and 8550 MW of total load \cite{Grigg:1999}. The pre-contingency
DC branch power fl{}ows have a mean of 113.8 MW, median of 109 MW,
and maximum of 396.1 MW. 

The second test system is a model of the 2004 winter peak Polish power
system that is available with MATPOWER \cite{Zimmerman:2011}. This
test system has 2896 branches (transmission lines and transformers),
2383 buses, and 24.6 GW of total load. The pre-contingency DC branch
power fl{}ows have a mean of 34 MW, median of 18.7 MW, and maximum
of 882.4 MW. For the Polish case, some of the transmission lines were
overloaded in the original system, so we increased line flow limits
to be the larger of the current limit and 1.05 times the maximum post-contingency
line flows in normal working condition for each line, after increasing
all loads by 10\%. This ensures that line limits are high enough to
serve all loads without load shedding after increasing the base case
demand by 10\%. 

Pre-contingency test cases were prepared for both systems using the
SCDCOPF method from Sec.~\ref{sec:Pre-contingency-Dispatch}, for
a range of load levels from 50\% to 119\%. 119\% was the highest load
factor for which SCDCOPF could find a solution without load shedding
in the RTS-96 case. For the Polish system, load could increase up
to 110\% without load shedding. However, we extended our study up
to 115\% for comparison, which caused less than 1\% load shedding
in the SCDCOPF solution for cases above 110\% load level. Finally,
Monte Carlo simulation was performed for both test networks, for each
load level.

\subsection{RTS-96 Results}

Fig.~\ref{fig:risk_RTS} shows cascading failure risk, in terms of
the expected value of blackout sizes, for two pre-contingency dispatch
conditions for the RTS case. Panel (a) shows the results after using
SCDCOPF at each load level. Each point on the graph shows the rolling
average of risk across 3 consecutive integer percentage load levels
(i.e., the datum at 90\% load is the average risk for 89\%, 90\%,
and 91\%). As previously mentioned, the risk associated with different
blackout sizes are separately presented. It is interesting to note
that small blackout risk is relatively uniform across all load levels,
with a peak at around 80\% load level, whereas large blackout risk
is largest at about 70\% load level, and decreases significantly as
the load level increases. Fig.~\ref{fig:risk_RTS} panel (b) shows
CF risk with a proportional dispatch method. To obtain these dispatch
cases, we took the 119\% load level case from SCDCOPF, and uniformly
decreased the loads and generators to each smaller load level. Interestingly,
this dispatch approach reduces risk substantially. Furthermore, what
risk remains is purely due to small blackouts (BO sizes < 5\%). The
pre-contingency dispatch in this case is obviously more expensive
that that from the SCDCOPF, which suggests that there is an important
tradeoff between generation dispatch costs and CF risk. 

In order to understand the reason behind this difference, we looked
at the power flows on five critical lines that connect the three areas
in the RTS-96 system. An outage on these lines can cause the system
to separate into islands. If this occurs and there is not enough generation
and time to allow generators to ramp up or down after the network
separates into islands, a large amount of load shedding may occur.
The line flow results in Table \ref{tab:Branch-flow-magnitude} show
that the flows are generally much higher at the 50\% and 75\% load
levels of the SCDCOPF dispatch than at the 119\% level. On the other
hand, for the proportional dispatch case, the power flows change more
uniformly as load changes. These results suggest that the SCDCOPF
algorithm is using more long-distance transmission at moderate load
levels, whereas at higher load levels important transmission corridors
are not loaded as close to their capacity.

\begin{table*}
\caption{\label{tab:Branch-flow-magnitude}Branch flow magnitude in five critical
lines connecting the three areas in RTS-96 with different load levels
and pre-contingency dispatch}

\begin{tabular}{cc|ccc|ccc}
\cline{3-8} 
 & \multicolumn{1}{c}{} & \multicolumn{6}{c}{Branch flow magnitude (MW)}\tabularnewline
\cline{3-8} 
 & \multicolumn{1}{c}{} & \multicolumn{3}{c|}{SCDCOPF dispatch} & \multicolumn{3}{c}{Proportional dispatch}\tabularnewline
\hline 
From (bus no.) & To (bus no.) & Load level: 50\% & Load level: 75\% & Load level: 119\% & Load level: 50\% & Load level: 75\% & Load level: 119\%\tabularnewline
\hline 
107 & 203 & 23.88 & 44.38  & 8.57 & 3.60  & 5.40  & 8.57\tabularnewline
\hline 
113 & 215 & 86.32 & 122.46 & 50.38 & 21.17 & 31.75  & 50.38\tabularnewline
\hline 
123 & 217 & 21.08 & 31.76 & 2.67 & 1.12  & 1.68  & 2.67\tabularnewline
\hline 
325 & 121 & 156.27 & 143.58 & 17.98 & 7.55  & 11.33  & 17.98\tabularnewline
\hline 
318 & 223 & \multicolumn{1}{c}{181.27} & 124.58 & 31.48 & 13.23 & 19.84 & 31.48\tabularnewline
\hline 
\end{tabular}
\end{table*}

\subsection{Polish System Results}

\begin{figure}
\begin{centering}
\includegraphics[width=0.9\columnwidth]{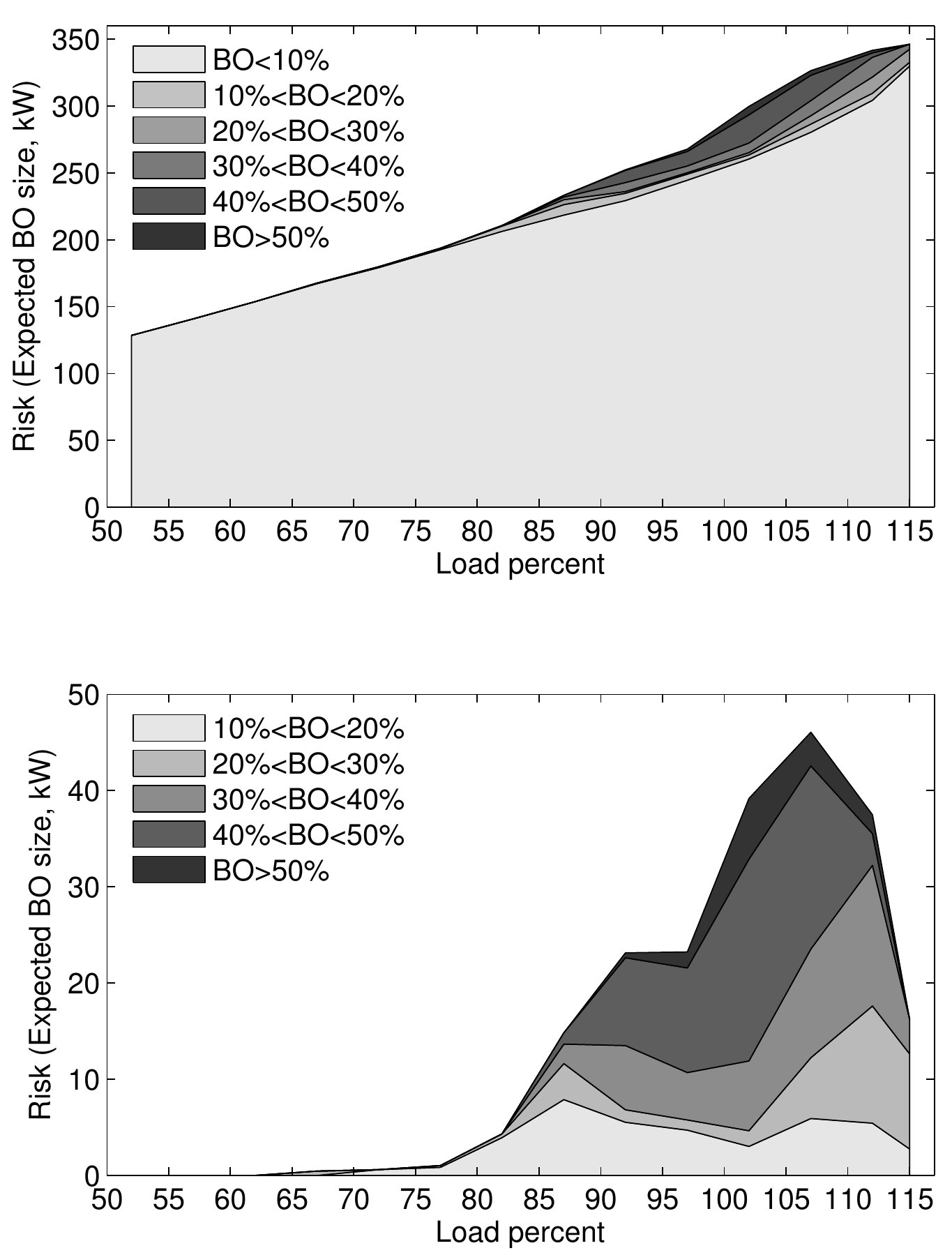}
\par\end{centering}

\caption{\label{fig:risk_polish}CF risk for the Polish grid; all blackouts
(top), big blackouts (bottom). }
\end{figure}

The same MC sampling approach was implemented on the Polish grid,
as explained in Section~\ref{sub:Test-Networks}. Here, we only use
SCDCOPF for the pre-contingency dispatch. Fig.~\ref{fig:risk_polish}
shows CF risk for all blackout sizes (top), and for only large blackouts
(bottom). This figure shows that there is a high risk associated with
small blackouts (BO sizes < 10\%), which increases uniformly with
load level. This result is largely due to the fact that the Polish
test case has numerous loads on radial lines, the failure of which
can cause load shedding in the down-stream system. Although larger
blackouts are less likely, their outcomes can be catastrophic. Therefore,
in this paper, we are mostly concerned with large blackouts from cascading
failure. Fig.~\ref{fig:risk_polish} (bottom) shows CF risk associated
with blackout sizes greater than 10\% separately for each 10\% interval.
We see that the pattern of these blackouts no longer changes uniformly
with load level. It is also worth noting that CF risk decreases for
load percentages higher than 110\% (the same cases for which some
load shedding, less than 1\%, occurs during the pre-contingency dispatch).

\section{Limitations of the Monte Carlo Approach and Future Work\label{sec:Monte-Carlo-Limitations}}

The expected value of a random variable (blackout size in our study)
can be found if we know the probability distribution function of that
random variable. Given our assumption that the results of a given
contingency are deterministic, if we can find all branch combinations
that cause a blackout, then CF risk is:

\begin{equation}
R(x)=\underset{c\in\mathcal{C}}{\sum}\textrm{Pr}(c)S(c,x)
\end{equation}
where $x$ is the initial system state, $c$ is a contingency, and
$\mathcal{C}$ is the set of all possible contingencies. $\textrm{Pr}(c)$
and $S(c,x)$ are, respectively, the probability and the blackout
size associated with each contingency $c$. In this paper, we computed
this index with Monte Carlo simulation, which required 30 hours of
computer time on a high performance platform \cite{VACC}, for each
load level. However, this amount of Monte Carlo simulation generated
primarily $n-2$ and a small number of $n-3$ contingencies. This
is due the fact that higher order outages are very low probability.
However, multiple contingencies do occur, and sometimes lead to large
cascading failures, adding to the importance of understanding risk
from these rare events. Finding even a small fraction of these rare,
but dangerous contingencies is a challenging task. Ref.~\cite{Eppstein:2012}
used Random Chemistry algorithm to identify a large collection of
$n-k$ contingencies that lead to cascading failure, but did not extend
the method to actually estimate CF risk. Implementing this extension
is a topic for future research.

\section{Conclusions\label{sec:Conclusions}}

In this paper, we study how cascading failure risk, measured using
the expected blackout size from Monte Carlo simulation, varies with
power system load level and pre-contingency dispatch. We used Monte
Carlo sampling of random branch outages to generate potential contingencies,
and a cascading failure simulator to evaluate the blackout sizes that
result. In order to understand the relative risk from different blackout
sizes, we separately measured risk associated with small, medium,
and large blackouts. The results indicate that, contrary to what one
might expect, risk does not necessarily increase with load level monotonically.
This is particularly true for larger blackouts (>5\% of load). The
results also show that the method used for pre-contingency generator
dispatch can play an important role in how risk changes with load
level. In one of our test systems, large blackout risk actually decreases
with load, after a threshold, when a security constrained optimal
power flow is used, whereas a more distributed method of dispatch
resulted in a risk profile that was both lower overall, and that increased
monotonically with load. While proportional dispatch resulted in lower
risk, it was also more expensive. This suggests that there is a tradeoff
between blackout risk and generator dispatch costs. 

Although Monte Carlo simulation is useful in cascading failure risk
estimation, because of the low probability, high impact nature of
cascading failure, the MC approach requires enormous computational
resources to obtain sufficiently low-variance risk estimates. Despite
their low probability, multiple-contingencies can trigger catastrophic
blackouts, and need to be studied carefully. Computing the risk of
cascading failure efficiently is the topic of future research.

\section*{Acknowledgment}

The authors gratefully acknowledge helpful conversations with Ian
Dobson, Daniel Kircshen and Steven Miller about these results, Maggie
Eppstein regarding our methods of risk estimation, and the Vermont
Advanced Computing Core (VACC), which is supported by NASA (NNX-08AO96G),
at the University of Vermont for providing High Performance Computing
resources. 

\bibliographystyle{IEEEtran}
\bibliography{Cascading_Failure}

\section*{Author Biographies}

\begin{IEEEbiography}[{\includegraphics[width=1in,height=1.25in,clip,keepaspectratio]{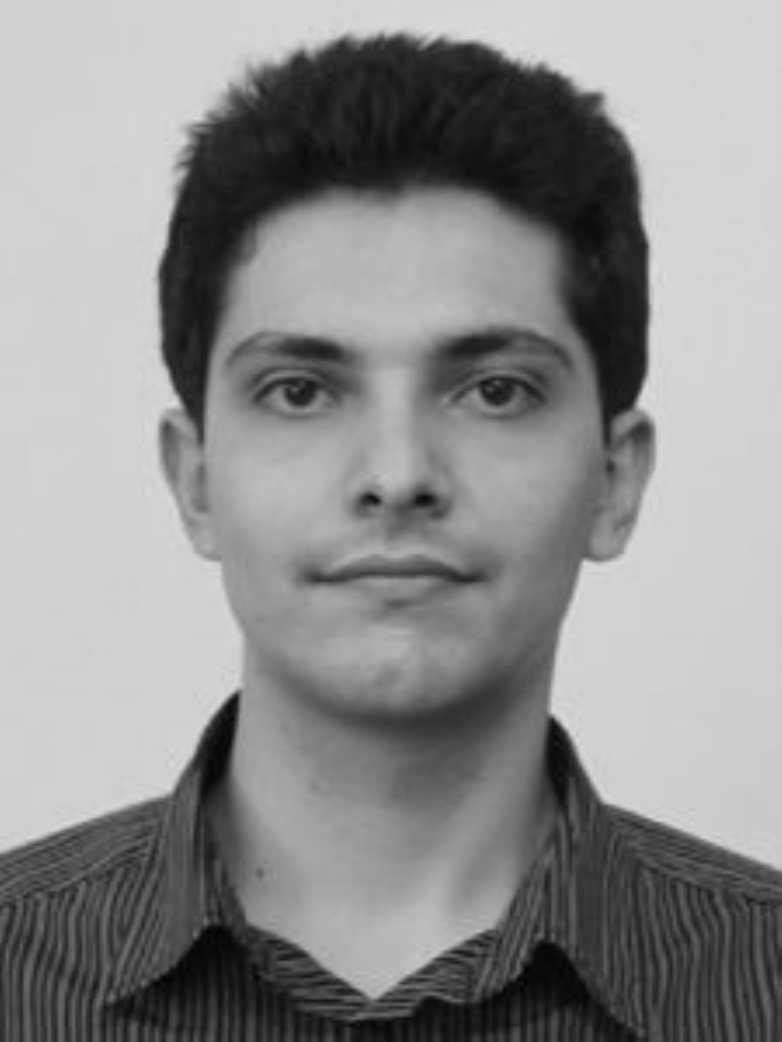}}]{Pooya Rezaei} (S'12) received the M.Sc. degree in Electrical Engineering from Sharif University of Technology, Tehran, Iran in 2010, and the B.Sc. degree in Electrical Engineering from University of Tehran, Tehran, Iran in 2008. Currently, he is pursuing the Ph.D. degree in Electrical Engineering at University of Vermont.  \end{IEEEbiography}

\begin{IEEEbiography}[{\includegraphics[width=1in,height=1.25in,clip,keepaspectratio]{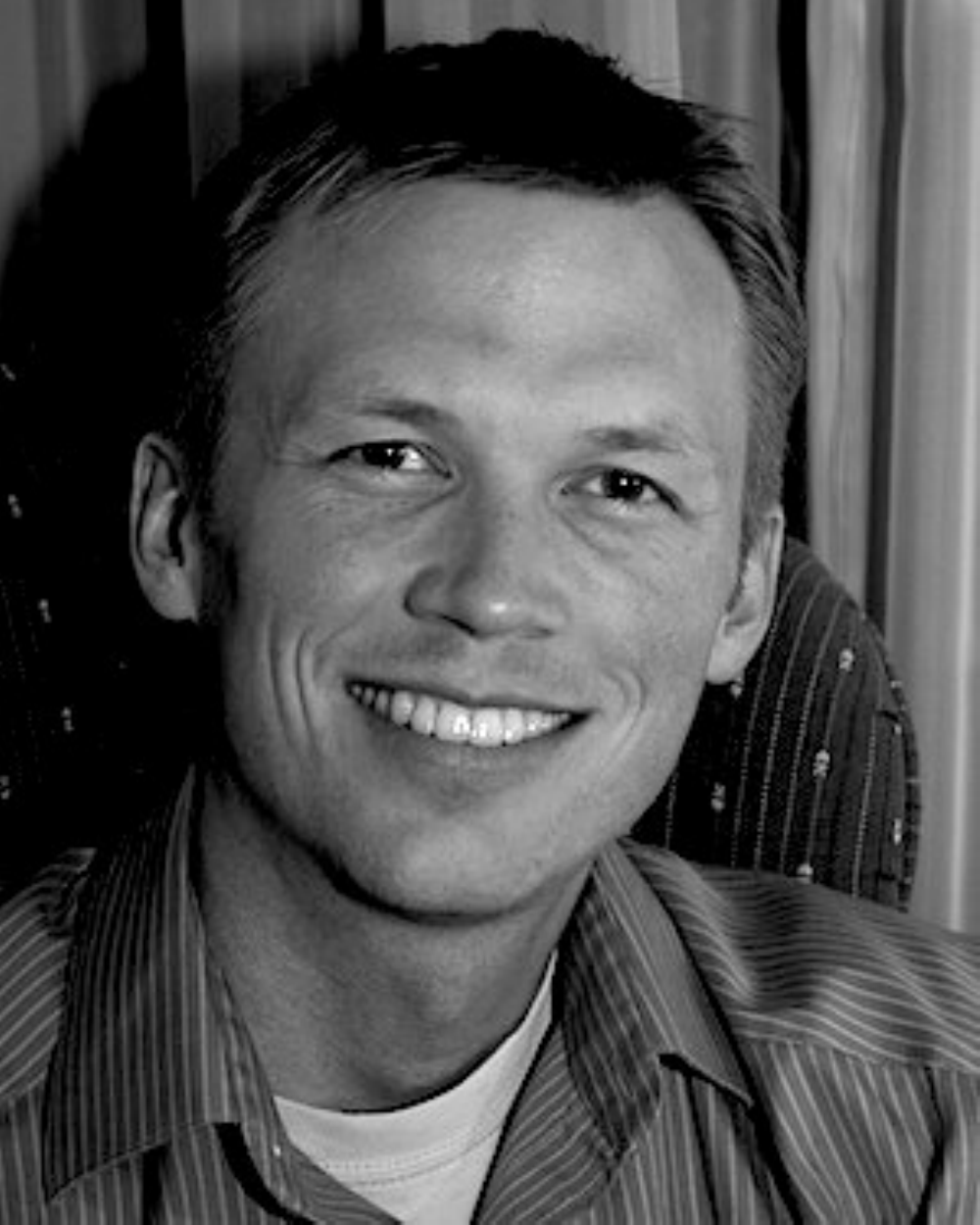}}]{Paul D. H. Hines} (S'96, M'07) received the Ph.D. degree in Engineering and Public Policy from Carnegie Mellon University in 2007, and the M.S. degree in Electrical Engineering from the University of Washington in 2001. He is currently an Assistant Professor in the School of Engineering at University of Vermont.
\end{IEEEbiography}
\vfill
\end{document}